\documentclass[12pt]{article}
\usepackage[top=3cm,bottom=3.5cm,left=3cm,right=3cm]{geometry}
\usepackage{graphicx}
\pdfoutput=1
\usepackage{mathrsfs}
\usepackage{amsfonts}
\usepackage{color}
\usepackage{amsmath}
\usepackage{setspace}
\usepackage{multirow}
\usepackage{cite}
\usepackage{ntheorem}
\usepackage{pstricks}
\usepackage{fancybox}
\usepackage{multirow}
\usepackage{indentfirst}
\usepackage{graphicx}
\usepackage{hyperref}

\usepackage{tikz}
\usetikzlibrary{arrows,positioning}

\theoremseparator{.}
\theorembodyfont{\normalfont}
\newtheorem{Assumption}{Assumption}
\newtheorem{Theorem}{Theorem}
\newtheorem{Definition}{Definition}

\newcommand\pubdate{\today}

\def\support{\footnote{Corresponding author: liux3771@umn.edu}}

\def\Title#1{\begin{center} {\Large #1 } \end{center}}
\def\Author#1{\begin{center}{  #1} \end{center}}
\def\Address#1{\begin{center}{ \it #1} \end{center}}

\newenvironment{Abstract}{\begin{quotation}  }{\end{quotation}}

\def\beq{\begin{equation}}
\def\eeq#1{\label{#1}\end{equation}}
\def\eeqn{\end{equation}}

\def\beqa{\begin{eqnarray}}
\def\eeqa#1{\label{#1}\end{eqnarray}}
\def\eeqan{\end{eqnarray}}

\let\bar=\overbar

\def\Dslash{\not{\hbox{\kern-4pt $D$}}}
\def\dslash{\not{\hbox{\kern-2pt $\del$}}}

\def\msb{{\bar{\ssstyle M \kern -1pt S}}}

\begin{document}
\begin{titlepage}
\pubdate
\vfill
\Title{Novel Criteria to Exclude the Surrogate Paradox and Their Optimalities}
\vfill
\Author{Yunjian Yin$^{a,b}$, Lan Liu$^{b,}$\support, Zhi Geng$^{a}$, and Peng Luo$^{c}$\\}
\Address{
$^{a}${\small School of Mathematical Sciences, Peking University, Beijing 100871, China}\\
$^{b}${\small School of Statistics, University of Minnesota, Minneapolis, Minnesota 55455, USA}
$^{c}${\small  Shenzhen University, Shenzhen, Guangdong 518060, China}}
\vfill
\begin{Abstract}
When the primary outcome is hard to collect, surrogate endpoint is typically used as a substitute. However, even when the treatment has a positive average causal effect (ACE) on the surrogate endpoint, which also has a positive ACE on the primary outcome, it is still possible that the treatment has a negative ACE on the primary outcome. Such a phenomenon is called the surrogate paradox and greatly challenges the use of surrogate. In this paper, we provide novel criteria to exclude the surrogate paradox. Unlike other conditions previously proposed, our conditions are testable since they only involve observed data. Furthermore, our criteria are optimal in the sense that they are sufficient and ``almost necessary" to exclude the paradox: if the conditions are satisfied, the surrogate paradox is guaranteed to be absent while if the conditions fail, there exists a data generating process with surrogate paradox that can generate the same observed data. That is, our criteria capture all the information in the observed data to exclude the surrogate paradox rather than relying on unverifiable distributional assumptions.\\
\noindent Key Words: Average causal effect; Optimality; Relative risk; Surrogate paradox.
\end{Abstract}
\vfill
\end{titlepage}
\def\thefootnote{\fnsymbol{footnote}}
\section{Introduction}
In many biomedical and econometric studies, the measurement of the primary endpoint may be expensive, inconvenient or infeasible to collect in a practical length of time. In such cases, surrogate variables or biomarkers are usually used as substitutes for the primary outcomes. For example, in cancer studies, the primary outcome is death, thus a surrogate endpoint is usually chosen to be tumor shrinkage or other laboratory measure to reduce the cost or the duration of the clinical trials \cite{fleming1996surrogate}.

A common misunderstanding of the surrogate endpoint is that a strong association between a surrogate and the primary outcome makes the surrogate a valid replacement for the primary outcome. However, such association does not justify the use of a surrogate: the effect of the treatment on the surrogate endpoint may not be a good predictor of the effect of the treatment on the outcome of interest. As illustrated in \cite{baker2003}, even with perfect correlation between the surrogate and the primary outcome, a positive treatment effect on the surrogate may still lead to a negative effect on the primary outcome.

Good surrogate is difficult to find \cite{burzykowski2006evaluation}, and the misuse of a surrogate may lead to severe consequences or even disasters \cite{fleming1996surrogate,moore1997deadly,manns2006surrogate}. For example, lipid levels, especially total cholesterol levels are significant predictors of cardiovascular-related mortality, thus the lipid lowering has been used as a surrogate for the reduction in mortality. However, Gordon \cite{gordon1995cholesterol} found that the use of cholesterol-lowering agents actually led to increase in overall mortality. Moore \cite{moore1997deadly} reported another example: anti-arrythmia drug Tamnbocor successfully suppresses arrythmia but resulted in the death of over 50,000 people.

To address this, different criteria have been proposed to evaluate the surrogate endpoints. The Prentice's criterion \cite{prentice1989surrogate}, also known as the statistical surrogate criterion, was the first operational criterion proposed to assess the validity of a surrogate. Given the surrogate, the criterion requires the conditional independence of the treatment and the primary endpoint. Intuitively, it requires the surrogate to capture all the relationship between the treatment and the primary outcome. Frangakis and Rubin \cite{frangakis2002principal} pointed out that the Prentice's criterion may not satisfy the property of causal necessity, i.e., the absence of the individual treatment effect on the surrogate indicates the absence of the treatment effect on the primary outcome. Instead, under the principal strata framework they proposed a principal surrogate criterion satisfies the property. Gilbert and Hudgens \cite{gilbert2008evaluating} introduced the concept of causal effect predictiveness (CEP) surface to evaluate a principal surrogate.  However, the principal surrogate criterion still lacks transportability in the sense that it relies on both the surrogate and the outcome to be available \cite{pearl2011principal}. Additionally, Pearl and Bareinboim \cite{pearl2010Transportability} presented counter-examples: a surrogate is a robust predictor of the treatment effect on the outcome but fails the principal surrogacy criterion as well as a surrogate satisfies the criterion but is useless as a predictor. Lauritzen \cite{lauritzen2004discussion} proposed a strong surrogate criterion using the language of graphical models. The key requirement of the strong surrogate criterion is that the treatment affects the primary outcome only through the surrogate endpoint. In other words, the strong surrogate fully mediates the effect of treatment on the outcome. By definition, a strong surrogate satisfies causal necessity and is also a principal surrogate.

However, all of the above criteria suffer from the surrogate paradox \cite{chen2007criteria}. As mentioned, a positive correlation between the surrogate and outcome does not make a good surrogate. The surrogate paradox implies that even a surrogate has a positive average causal effect (ACE) on the true endpoint, a positive average treatment effect on the surrogate may still lead to a negative treatment effect on the primary outcome. Chen et al. \cite{chen2007criteria} also demonstrated that the surrogate paradox can manifest even in randomized trials: although the treatment is randomized, the surrogate and primary outcome may still be confounded by some unmeasured confounders. To make progress, they provided sufficient conditions on the joint distribution of unmeasured confounders, surrogate and primary outcome to avoid the paradox. However, the conditions they proposed involve the unobserved confounders thus are untestable. Ju and Geng \cite{ju2010criteria} and Wu et al. \cite{wu2011sufficient} considered using the distributional causal effects (DCE) as the measure of causal effects and derived sufficient conditions to avoid the paradox.
Their conditions involve the outcome distributions of both the treatment and control arms, which are generally not completely observed in a case that involves surrogate endpoints.

In this paper, we propose novel criteria to avoid the surrogate paradox in randomized controlled trials with only surrogate being collected. We utilize the outcome information of the control arm from a previous study. We show that it is not enough to avoid the surrogate paradox merely with the ACE of surrogate on outcome being positive, instead, we require its magnitude to pass certain positive threshold.

Our work of avoiding surrogate paradox is novel in several aspects. First, our data setting is very common and useful in practice. When it is of interest to assess the effect on the primary outcome of a new treatment versus a reference treatment, our conditions enable us to predict the sign of the treatment effect on the primary outcome in the new randomized controlled trial without collection of the outcome. Second, our conditions for the surrogate has transportability in the sense that the information about the causal relationship between the surrogate and the primary outcome from a previous study could be used in the current study. Third, our methods could be used to exclude the existence of surrogate paradox using the observed data rather than relying on unverifiable distributional assumptions. Finally, our criteria are optimal in the sense that they are sufficient and ``almost necessary" if the conditions are satisfied, the surrogate paradox is guaranteed to be absent while if the conditions fail, there exists a data generating process (a set of parameters) with surrogate paradox that can generate the same observed data. That is, our criteria capture all the information in the observed data to exclude the surrogate paradox.

The paper is organized as follows. In Section 2, we introduce the motivating example. The notation and assumptions are introduced in Section 3. The surrogate paradox is reviewed in Section 4. In Section 5, we derive the sharp bounds for the ACE of the treatment on the primary outcome within the class of strong surrogate endpoints. Based on that, a novel condition is proposed to avoid the surrogate paradox. In Section 6, we extend the condition to non-strong surrogate endpoint and also extend to the situation where the causal effect is evaluated using causal relative risk (CRR) scale. We apply our method in the motivating example in Section 7. The paper concludes with a discussion in Section 8.

\section{Motivating Example}
Hypertension is an important risk factor and a major cause for cardiovascular diseases such as stroke, microvascular disease, myocardial infarction and for death\footnote{For example, Lawes et al. \cite{lawes2008global} found that about 54\% of stroke, 47\% of ischaemic heart disease, 75\% of hypertensive disease, and 25\% of other cardiovascular disease worldwide were attributable to high blood pressure.}. In total, about 7.6 million (13.5\%) of all deaths and 92 million (6.0\%) of all disability adjusted life years (DALYs) in the year 2001 were due to high blood pressure as a cause of these diseases \cite{lawes2008global}. Hypertension is even more common in people with type 2 diabetes than in the general population. Stamler et al. \cite{stamler1993diabetes} suggested that diabetes mellitus increases the risk of cardiovascular disease by a factor of two to three at every level of systolic blood pressure (SBP).

The collection of the development of cardiovascular disease or death usually require long duration of the study. Thus in many countries may approve anti-hypertension drug based on the evidence of surrogate efficacy; that is, such drugs reduce blood pressure \cite{fleming1996surrogate}. However, based on a population-based case-control study, Psaty et al. \cite{psaty1995risk} suggested that the anti-hypertension drug calcium channel blockers may be associated with an increased risk for myocardial infarction among hypertensive patients. Such undesirable association suggests the need for careful investigation of  the use of blood pressure as a valid surrogate.

In order to promote the benefit of blood pressure lowering and blood pressure management, the Seventh Report of the Joint National Committee on Prevention, Detection, Evaluation, and Treatment of High Blood Pressure (JNC 7) recommended starting drug treatment in patients with diabetes who have SBPs of 130 mm Hg or higher, with a treatment goal of reducing SBP to below 130 mm Hg.

To justify such recommendations, the Action to Control Cardiovascular Risk in Diabetes (ACCORD) group conducted a randomized trial at 77 clinical sites from seven networks in the United States and Canada \cite{accord2010effects}. The ACCORD blood pressure trial aimed to compare an intensive therapy (a therapeutic strategy that targets a SBP $<120$ mmHg) with a standard therapy (a strategy that targets a SBP of $<140$ mmHg). A total of 4733 participants with type 2 diabetes were randomly assigned to intensive or standard therapy and were followed for 4.7 years on average. Baseline characteristics were generally similar between the two groups. After 1 year, the mean SBP was 119.3 mm Hg in the intensive therapy group and 133.5 mm Hg in the standard therapy group. Since the long term benefit of blood pressure lowering has been used as evidence for policy making and clinical recommendation \cite{macmahon1994blood,ogden2000long}, the primary outcome in our analysis is the death from cardiovascular causes over 10 years.

Although the  intensive therapy was new, the standard strategy has been adopted in previous trials. The United Kingdom Prospective Diabetes Study (UKPDS) is a prospective observational study \cite{adler2000association}. The study population consists of 4801 hypertensive participants in an epidemiological component and 1148 patients with hypertension from a blood control component. The patients from the UKPDS blood control component received treatment to control their blood pressure with a $\beta$ blocker or an angiotensin converting enzyme inhibitor to control their SBP to $<$150 mmHg or to even less tight control. If the target blood pressure was not met, additional agents were prescribed including a loop diuretic, a calcium channel blocker, and a vasodilator. The median follow-up time for all cause mortality was 10.5 years.

Thus, available information comes from three resources: (1) the data of both surrogate endpoint and primary endpoint of patients using standard therapy from the UKPDS cohort (2) the data of surrogate endpoints of patients using standard or intensive therapy from the ACCORD randomized controlled trial and (3) the causal relationship between elevated blood pressure (surrogate endpoint) and the mortality (primary outcome).

\section{Notation and Assumptions}
Let $T$ denote the treatment variable, $Y$ the primary outcome and $S$ the surrogate endpoint. We restrict $T$, $S$ and $Y$ to be binary throughout the paper. Let $X$ denote the observed covariates and let $U$ denote the unobserved confounder that affects both $S$ and $Y$. We do not impose any restrictions on the unmeasured confounder $U$ and it can be discrete or continuous. Let $t$ denote the possible value $T$ could take ($t = 1$ for treatment and $t = 0$ for placebo). We assume larger values indicate better results for both $S$ and $Y$. Under the Stable Unite Treatment Value Assumption (SUTVA) \cite{rubin1980randomization}, let $Y_{ts}$ denote the potential outcomes of the primary endpoint if the treatment and the surrogate were set to $T =t$ and $S=s$ by an external intervention on $T$ and $S$. Let $S_{t}$ denote the potential outcome of surrogate if the treatment was set to $T=t$. We may also use the notation $Y_{T=t}$ as the potential primary outcome when the intervention is only to set $T=t$. Let $\mathrm{ACE}(T\rightarrow Y)$ denote the ACE of the treatment $T$ on the outcome $Y$, i.e., $\mathrm{ACE}(T\rightarrow Y)=P(Y_{T=1}=1)-P(Y_{T=0}=1)$. Likewise, we define $\mathrm{ACE}(T\rightarrow S)$ and $\mathrm{ACE}(S\rightarrow Y)$ as the ACE of $T$ on $S$ and $S$ on $Y$ respectively. Similarly, we can define the stratified ACE based on the values of observed covariates $X$. For example, define $\mathrm{ACE}(T\rightarrow Y|X=x)$ to be the ACE of treatment $T$ on outcome $Y$ among the subgroups of people whose covariates $X$ take value $x$. We make the following randomization assumption throughout the paper.

\begin{Assumption}\label{A1}
(Randomization) $T\bot(Y_{00},Y_{01},Y_{10},Y_{11},S_{0},S_{1},U)$.
\end{Assumption}

The randomization assumption can be relaxed by the conditional ignorability assumption: suppose $T$ is independent of all potential outcomes and unmeasured confounder $U$ conditional on $X$, i.e., $T\bot(Y_{00},Y_{01},Y_{10},Y_{11},S_{0},S_{1},U)|X$.

To ease the illustration, we first demonstrate our methods within the class of strong surrogates.

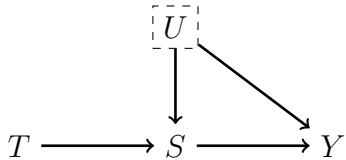
\begin{figure}
\centering
\begin{tikzpicture}
\node[text centered] (t) {$T$};
\node[right = 1.5 of t, text centered] (s) {$S$};
\node[right=1.5 of s, text centered] (y) {$Y$};
\node[draw, rectangle, dashed, above = 1 of s, text centered] (u) {$U$};

\draw[->, line width= 1] (t) --  (s);
\draw [->, line width= 1] (s) -- (y);
\draw[->,line width= 1] (u) --(s);
\draw[->,line width= 1] (u) -- (y);
\end{tikzpicture}
\medskip
\caption{Causal diagram of the strong surrogate $S$ for the effect of the treatment $T$ on outcome $Y$.\label{DAG}}
\end{figure}

\begin{Assumption}\label{A2}
(Strong Surrogate) The causal diagram in Figure \ref{DAG} is valid.
\end{Assumption}

A surrogate is a strong surrogate for the effect of the treatment $T$ on the outcome $Y$ if the entire treatment effect on the outcome goes through the surrogate \cite{lauritzen2004discussion}. In the mediation literature, $S$ is called to fully mediate the effect from $T$ to $Y$. We can also interpret $T$ as an instrumental variable for surrogate $S$ on $Y$.

Under assumption \ref{A2}, it is easy to see that $Y_{ts}$ does not depend on $t$, i.e., $Y_{ts}=Y_{s}$. We will omit the subscript $t$ and denote the potential outcome as $Y_{s}$ hereafter, unless otherwise specified. Hence there are four possible values for the vector $(Y_{0},Y_{1})$ with $Y_{s}\in\{0,1\}$ for $s=0,1$. Additionally, there are also four possible values for the vector $(S_{0},S_{1})$ with $S_{t}\in\{0,1\}$ for $t=0,1$. Thus, we have $4\times 4=16$ possible values concerning the vector of the potential outcomes $(Y_{0},Y_{1},S_{0},S_{1})$. Let $q_{ij}\geq0$ denote the proportion of getting each possible value of $(Y_{0},Y_{1},S_{0},S_{1})$ in the whole population for $i,j=0,\ldots,3$ (as shown in Table \ref{tb: types}), e.g., $q_{01}=P(Y_{S=0}=0,Y_{S=1}=1,S_{T=0}=0,S_{T=1}=0)$.
\renewcommand\arraystretch{1.6}
\begin{table}
\caption{Probabilities of sixteen potential outcomes types in the strong surrogate scenario\label{tb: types}}
\centering
\begin{tabular}{|l|l|l|l|l|}
  \hline
   & $Y_{S=0}=0$ & $Y_{S=0}=0$ & $Y_{S=0}=1$ & $Y_{S=0}=1$ \\
   & $Y_{S=1}=0$ & $Y_{S=1}=1$ & $Y_{S=1}=0$ & $Y_{S=1}=1$\\
    \hline
  $S_{T=0}=0,S_{T=1}=0$ & $q_{00}$ & $q_{01}$ & $q_{02}$ & $q_{03}$ \\
   \hline
  $S_{T=0}=0,S_{T=1}=1$ & $q_{10}$ & $q_{11}$ & $q_{12}$ & $q_{13}$ \\
   \hline
  $S_{T=0}=1,S_{T=1}=0$ & $q_{20}$ & $q_{21}$ & $q_{22}$ & $q_{23}$ \\
   \hline
  $S_{T=0}=1,S_{T=1}=1$ & $q_{30}$ & $q_{31}$ & $q_{32}$ & $q_{33}$ \\
  \hline
\end{tabular}
\end{table}
\renewcommand\arraystretch{0.8}

Thus, the ACE of $T$ on $Y$ can be expressed as:
\begin{eqnarray*}\label{ACETY}
&&\mathrm{ACE}(T\rightarrow Y)\notag \\
&=&P(Y_{T=1}=1)-P(Y_{T=0}=1) \notag \\
    &=&\sum_{s=0,1}P(Y_{T=1}=1,S_{1}=s)-\sum_{s=0,1}P(Y_{T=0}=1,S_{0}=s) \notag \\
    &=& P(Y_{S=0}=1,S_{1}=0)+P(Y_{S=1}=1,S_{1}=1)-P(Y_{S=0}=1,S_{0}=0)-P(Y_{S=1}=1,S_{0}=1) \notag \\
    &=& q_{22}+q_{11}-q_{12}-q_{21}.
\end{eqnarray*}
\noindent Generally, in a randomization study where a surrogate endpoint $S$ is collected as opposed to the outcome $Y$, the $\mathrm{ACE}(T\rightarrow Y)$ can not be estimated without strong modeling assumptions. The surrogate paradox reviewed in the next section indicates that even the sign of $\mathrm{ACE}(T\rightarrow Y)$ is not easy to predict in general.

\section{Surrogate Paradox}\label{Surrogate Paradox}
Chen et al. \cite{chen2007criteria} defined the phenomenon of surrogate paradox as follows:
\begin{itemize}
  \item [(a)] A treatment has a positive ACE on a surrogate which in turn has a positive ACE on a primary outcome, but the treatment has a negative ACE on the primary outcome, i.e., $\mathrm{ACE}(T\rightarrow S)>0$ and $\mathrm{ACE}(S\rightarrow Y)>0$, but $\mathrm{ACE}(T\rightarrow Y)<0$.
\end{itemize}

We can also have the following surrogate paradox phenomenon by redefining $T^{*}=1-T$, $S^{*}=1-S$, $Y^{*}=1-Y$:
\begin{itemize}
  \item [(b)] A treatment has a positive ACE on a surrogate which in turn has a negative ACE on a primary outcome, but the treatment has a positive ACE on the primary outcome, i.e., $\mathrm{ACE}(T\rightarrow S)>0$ and $\mathrm{ACE}(S\rightarrow Y)<0$, but $\mathrm{ACE}(T\rightarrow Y)>0$;
  \item [(c)] A treatment has a negative ACE on a surrogate which in turn has a positive ACE on a primary outcome, but the treatment has a positive ACE on the primary outcome, i.e., $\mathrm{ACE}(T\rightarrow S)<0$ and $\mathrm{ACE}(S\rightarrow Y)>0$, but $\mathrm{ACE}(T\rightarrow Y)>0$;
  \item [(d)] A treatment has a negative ACE on a surrogate which in turn has a negative ACE on a primary outcome, but the treatment has a negative ACE on the primary outcome, i.e., $\mathrm{ACE}(T\rightarrow S)<0$ and $\mathrm{ACE}(S\rightarrow Y)<0$, but $\mathrm{ACE}(T\rightarrow Y)<0$.
\end{itemize}

Without loss of generality, we focus on the definition (a) of the surrogate paradox and hereafter we assume that the $\mathrm{ACE}(T\rightarrow S)$ and $\mathrm{ACE}(S\rightarrow Y)$ are both positive.

For a good surrogate, the sign of $\mathrm{ACE}(T\rightarrow Y)$ should be the same as the sign of the product of $\mathrm{ACE}(T\rightarrow S)$ and $\mathrm{ACE}(S\rightarrow Y)$. But the phenomena above violates this.
For illustration, Chen et al. \cite{chen2007criteria} considered the example of antiarrhythmic drugs. They assumed that the surrogate is a strong surrogate and that all the variables including the confounder $U$ are binary. The probabilities they considered are $P(T=1)=0.5$, $P(U=1)=0.7$ and others given in Table \ref{tb: Chen_eg}. From these probabilities we can have $\mathrm{ACE}(T\rightarrow S)=0.6220$, $\mathrm{ACE}(S\rightarrow Y)=0.3010$, $\mathrm{ACE}(T\rightarrow Y)=-0.0491$. This manifests the surrogate paradox (a).

\renewcommand\arraystretch{2}
\begin{table}
\caption{Example of antiarrhythmic drugs in \cite{chen2007criteria}\label{tb: Chen_eg}}
\centering
\begin{tabular}{|cccccc|}
  \hline
   & \multicolumn{2}{c}{$P(S=1|U,T)$}  && \multicolumn{2}{c|}{$P(Y=1|U,S)$}   \\
\cline{2-3}\cline{5-6}
    & $T=0$ & $T=1$ && $S=0$ & $S=1$ \\
    \hline
  $U=0$ & 0.98 & 0.79 && 0.00 & 0.98 \\
  $U=1$ & 0.02 & 0.99 && 0.98 & 0.99 \\
  \hline
\end{tabular}
\end{table}
\renewcommand\arraystretch{0.8}

The fundamental reason for the surrogate paradox is the absence of the primary outcome $Y$ and the presence of the unmeasured confounder $U$. Since the primary outcome is not fully observed, the presence of the latent factor that affects both the surrogate and the outcome would lead to completely unexpected results. It can be proved that, when $U$ is absent, we have $\mathrm{ACE}(T\rightarrow Y)=\mathrm{ACE}(T\rightarrow S)\mathrm{ACE}(S\rightarrow Y)$ for a strong surrogate $S$. This indicates that without the unmeasured confounder $U$, not only the sign, but also the magnitude of $\mathrm{ACE}(T\rightarrow Y)$ can be predicted.

As demonstrated in \cite{baker2003}, a strong positive correlation between the surrogate and the outcome does not make a good surrogate. Even if replacing the association with a causal effect, the mere knowledge of the sign of $\mathrm{ACE}(S\rightarrow Y)$ cannot exclude the surrogate paradox. 
In the next section, we will show that if the magnitude of $\mathrm{ACE}(S\rightarrow Y)$ is bigger than certain threshold, the existence of the surrogate paradox could be excluded using the observed data. We will also show that such criterion to exclude the surrogate paradox satisfies certain optimality.

\section{Excluding the Surrogate Paradox}\label{boundsection}
\subsection{Deriving the bounds with linear programming}\label{subsec: bounds}
In order to provide condition to exclude the surrogate paradox, we first derive the sharp bounds of $\mathrm{ACE}(T\rightarrow Y)$ with the available data $\big\{P(S,Y|T=0),P(S|T=1),\mathrm{ACE}(S\rightarrow Y)\big\}$. To simplify notation, let $\gamma=\mathrm{ACE}(S\rightarrow Y)$ denote the ACE of the surrogate $S$ on the outcome $Y$. Under the assumptions \ref{A1} and \ref{A2}, we have the following restrictions on the distributional parameters of the potential outcomes:
\begin{equation}\label{constraint1}
\left\{
  \begin{array}{ll}
    P(Y=0,S=0|T=0)=q_{00}+q_{01}+q_{10}+q_{11},  \\
    P(Y=1,S=0|T=0)=q_{02}+q_{03}+q_{12}+q_{13},  \\
    P(Y=0,S=1|T=0)=q_{20}+q_{22}+q_{30}+q_{32},  \\
    P(S=0|T=1)=q_{00}+q_{01}+q_{20}+q_{21}+q_{02}+q_{03}+q_{22}+q_{23},  \\
    \gamma=q_{01}+q_{11}+q_{21}+q_{31}-q_{02}-q_{12}-q_{22}-q_{32}, \\
    \sum_{i=0}^{3}\sum_{j=0}^{3}q_{ij}=1,
  \end{array}
\right.
\end{equation}
By the linear programming method, we have following bounds for the $\mathrm{ACE}(T\rightarrow Y)$.

\begin{Theorem}\label{Th1}
Under the assumptions \ref{A1} and \ref{A2}, if $\big\{P(S,Y|T=0),P(S|T=1),\gamma\big\}$ is known, we can have the sharp bounds for $\mathrm{ACE}(T\rightarrow Y)$, denote as $[L,U]$, where
$$L =\max\left(
                               \begin{array}{c}
                                 L_{1} \\
                                 L_{2} \\
                                 L_{3} \\
                                 L_{4} \\
                                 L_{5} \\
                                 L_{6} \\
                                 L_{7} \\
                               \end{array}
                             \right)\equiv
\max\left(
               \begin{array}{c}
                 -P(Y=1,S=0|T=0)-P(S=0|T=1) \\
                 -P(Y=1|T=0) \\
                 -P(Y=1,S=1|T=0)-P(S=1|T=1) \\
                 -\gamma-P(S=1|T=1)-P(Y=1,S=0|T=0)\\
                 -\gamma-P(Y=0,S=1|T=0)-2P(Y=1,S=0|T=0)\\
                 \gamma-2P(Y=1,S=1|T=0)-P(Y=0,S=0|T=0)\\
                 \gamma-P(S=0|T=1)-P(Y=1,S=1|T=0)\\
               \end{array}
             \right),$$
and
$$U=\min\left(
                               \begin{array}{c}
                                 U_{1} \\
                                 U_{2} \\
                                 U_{3} \\
                                 U_{4} \\
                                 U_{5} \\
                                 U_{6} \\
                                 U_{7} \\
                               \end{array}
                             \right)\equiv \min\left(
               \begin{array}{c}
                 P(Y=0,S=0|T=0)+P(S=0|T=1) \\
                 P(Y=0|T=0) \\
                 P(Y=0,S=1|T=0)+P(S=1|T=1) \\
                 \gamma+P(Y=1,S=0|T=0)+2P(Y=0,S=1|T=0)\\
                 \gamma+P(Y=0,S=1|T=0)+P(S=0|T=1)\\
                 -\gamma+2P(Y=0,S=0|T=0)+P(Y=1,S=1|T=0)\\
                 -\gamma+P(Y=0,S=0|T=0)+P(S=1|T=1)\\
               \end{array}
             \right).
$$
\end{Theorem}
Note that the bounds in Theorem \ref{Th1} are sharp in the sense that they can be obtained by certain choice of joint distribution of full data $(T,S,Y,U)$. Both the upper bound and the lower bound only involve the observed data distribution and  can be estimated using nonparametric method. As seen from the results, the width of bounds does not have a fixed value and how informative the bounds depend on the observed data. Specifically, when $L_{7}=U_{7}$, i.e. $\gamma=\{1+P(Y=0,S=0|T=0)+P(Y=1,S=1|T=0)\}/2$, the bounds vanish to a single point. In such case, the causal effect of $T$ on $Y$ can be identified.

\subsection{Exclude the surrogate paradox using the observed data}

In the absence of the paradox, when $\mathrm{ACE}(T\rightarrow S)>0$ and $\mathrm{ACE}(S\rightarrow Y)>0$, we have $\mathrm{ACE}(T\rightarrow Y)>0$. Since $\mathrm{ACE}(T\rightarrow Y)\geq L$, the condition that $L>0$ excludes the surrogate paradox. Note that $\max\{L_{1},L_{2},L_{3},L_{4},L_{5}\}\leq 0$, so $L>0$ is equivalent to $\max\{L_{6},L_{7}\}>0$. That is also equivalent to
\begin{equation}\label{SN}
\gamma>\min\left(
               \begin{array}{c}
                2P(Y=1,S=1|T=0)+P(Y=0,S=0|T=0)\\
                P(S=0|T=1)+P(Y=1,S=1|T=0)\\
               \end{array}
             \right).
\end{equation}
Thus, \eqref{SN} can be used as a criterion to exclude the paradox, that is, when \eqref{SN} is satisfied, the surrogate paradox does not exist.

\begin{Theorem}\label{Th2}
Under the assumptions \ref{A1} and \ref{A2}, if $\big\{P(S,Y|T=0),P(S|T=1),\gamma\big\}$ is known, the surrogate paradox can be excluded if (\ref{SN}) holds. 
\end{Theorem}

Note that the right hand side of \eqref{SN} is testable since it only involves the available data. This is an appealing feature of the criterion. In contrast, if the exclusion of the surrogate paradox completely relies on unverifiable conditions, there is little support from the data to justify the use of the surrogate as a good replacement for the outcome.
An extreme example of such unverifiable condition would be assuming $\mathrm{ACE}(T\rightarrow Y)>0$, which is the same as suggesting the use of surrogate without justification.  
\begin{figure}
\includegraphics[width=0.8\textwidth]{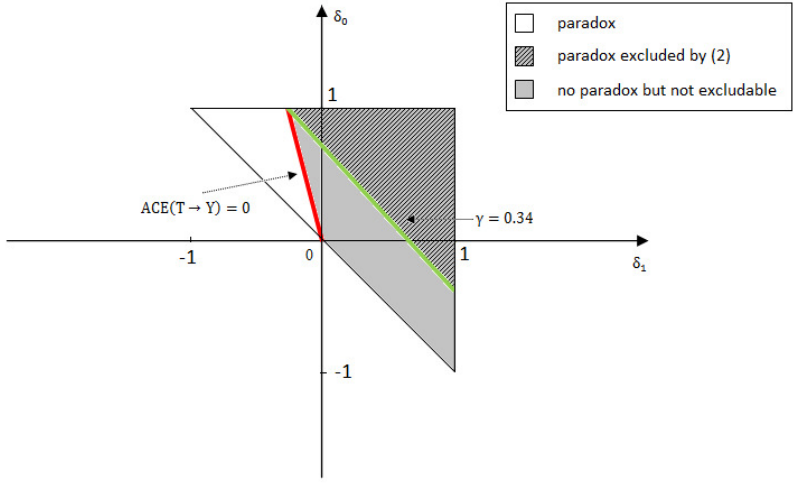}
  \caption{Partition of the parameter space of $(\delta_0,\delta_1)$}\label{fig: projection}
\end{figure}

To further illustrate the exclusion of surrogate paradox by \eqref{SN}, we demonstrate in Figure \ref{fig: projection} the partition of the parameter space based on the manifestation of the paradox. For the ease of illustration, we reduce the dimension of parameters by setting $P(U=1)=0.5$, $P(S=1|T=1,U=0)-P(S=1|T=0,U=0)=0.7$ and $P(S=1|T=1,U=1)-P(S=1|T=0,U=1)=0.2$. Let $\delta_u=P(Y=1|S=1,U=u)-P(Y=1|S=0,U=u)$ for $u=0,1$. Figure \ref{fig: projection} depicts the partition of parameter space of $\delta_0$ and $\delta_1$ according to whether the paradox manifests. Note that the outer triangle, bounded by the horizontal and vertical lines at 1 and a line with slope -1 through the origin, denotes the parameter space for $\delta_0$ and $\delta_1$ such that $\mathrm{ACE}(T\rightarrow S)>0$ and $\mathrm{ACE}(S\rightarrow Y)>0$. Thus, any point in this area corresponds to a set of full data distribution of $(U,T,S,Y)$ with $\mathrm{ACE}(T\rightarrow S)>0$ and $\mathrm{ACE}(S\rightarrow Y)>0$. The white triangle area denotes the parameter space where the surrogate paradox manifests while the grey area (including both the shaded and non shaded area) denotes the parameter space where the surrogate paradox is absent. Since the distribution of the full data $(U,T,S,Y)$ is not fully identifiable, one can not pin point where the parameters are in the figure based on observed data. The green line denotes the contour line of $\gamma=0.34$ and the shaded area denotes the parameter space such that $\gamma>0.34$, where $0.34$ is the threshold of $\gamma$ specified by of the right hand side of \eqref{SN}. Thus by theorem \ref{Th2}, if the parameters fall into the shaded area, the surrogate paradox can be excluded by the observed data alone.


Apart from the testability, our criterion also has the following optimality.

\begin{Definition}\label{def: optimal}
A criterion to exclude the surrogate paradox is {\it optimal} if (i) when the criterion is satisfied, the surrogate paradox is absent (ii) when the criterion is not satisfied, one can always find a data generating mechanism that yields the same observed data but suffers from the surrogate paradox. That is, one cannot exclude the possibility of surrogate paradox according to the observed data.
\end{Definition}

Intuitively, an ideal criterion to exclude the surrogate paradox will be based on a sufficient and ``almost necessary" condition. The sufficiency gives the condition enough strength to rule out surrogate paradox: if the condition is satisfied, the surrogate paradox is guaranteed to be absent. The ``almost necessity" gives the condition enough sharpness to hold as long as the observed data could rule out the possibility of surrogate paradox: if the condition fails, there exists a data generating process (a set of parameters) with surrogate paradox that can generate the same observed data.

We have demonstrated the sufficiency of our condition in Figure \ref{fig: projection}. For a graphical illustration of ``almost necessity", note that for any full data distribution corresponds to a point in the grey unshaded area, there is a full data distribution corresponds to a point in the white triangle area that can generate the same observed data. Thus, if the parameters fall into the grey unshaded area, the surrogate paradox cannot be excluded based on the observed data.


The ``almost necessity" differs from necessity in the sense that a necessary condition would require a criteria to rule out the possibility of surrogate paradox whenever it is absent. Such necessity is impossible to achieve since the underlying data generating mechanism is generally not identifiable. More specifically, we can only identify a set of data-generating process that is consistent with the observed data. If and only if none of these data generating mechanisms has surrogate paradox, the criterion enable us to exclude surrogate paradox. That said, the optimality in the definition \ref{def: optimal} requires a criterion to capture all the information in the {\it observed data} to exclude the surrogate paradox.

The optimality of condition \eqref{SN} is given in the following theorem.  Additionally, in the appendix, we show that all the previously proposed criteria to exclude the surrogate paradox are not optimal.

\begin{Theorem}\label{Th3}
Under the assumptions \ref{A1} and \ref{A2}, \eqref{SN} is an optimal criterion to exclude the surrogate paradox.
\end{Theorem}

\vspace{-2mm}
The information of $\gamma$ plays an important role in both theorems \ref{Th2} and \ref{Th3}. However, such information is not restricted to the magnitude of $\gamma$, but can also be a bound for $\gamma$. More specifically, in the appendix, we explain how to exclude paradox if we know a range of $\gamma$, e.g., $a<\gamma< b$ for some constants $a$, $b$. Also, we show that the paradox cannot be excluded with the sign of $\gamma$ being positive, i.e., $a=0$ and $b=+\infty$, regardless of what observed data one obtain. Thus, it is not enough to avoid the surrogate paradox merely with the ACE of surrogate on outcome being positive, instead, we require its magnitude to pass certain positive threshold.


Note that our bounds and conditions to avoid the surrogate paradox are given for the overall ACE instead of the stratified ACE. All the bounds and conditions we proposed can be easily extended to the cases where the covariates are discrete and of low dimensions. When the covariates $X$ is continuous or high dimensional, a parametric model for the ACE is needed due to the curse of dimensionality.

\subsection{Some illustrative examples}
We illustrate the Theorem \ref{SN} with three examples in this subsection. For simplicity, the confounder $U$ is chosen to be binary for each example.

\begin{description}

\vspace{-2mm}  \item[Example 1: (sufficiency)] We first revisit the example in \cite{chen2007criteria} which was also described in Section \ref{Surrogate Paradox}. With all the information about the joint distribution of $(T,S,U,Y)$, we have $\mathrm{ACE}(T\rightarrow S)=0.6220$, $\gamma=\mathrm{ACE}(S\rightarrow Y)=0.3010$, $\mathrm{ACE}(T\rightarrow Y)=-0.0491$, and the right hand side of (\ref{SN}) is 0.3720. Thus the inequality (\ref{SN}) does not hold and the surrogate paradox occurs.

\vspace{-2mm}
\item[Example 2: (sufficiency)] Let $P(T=1)=0.5$, $P(U=1)=0.3$ and the conditional probabilities are given in Table \ref{Ex3}. Then we have $\mathrm{ACE}(T\rightarrow S)=0.4420$, $\gamma=\mathrm{ACE}(S\rightarrow Y)=0.5750$, $\mathrm{ACE}(T\rightarrow Y)=0.2726$, and the right hand side in (\ref{SN}) is 0.4864. Thus the inequality (\ref{SN}) holds and the surrogate paradox is avoided.
\renewcommand\arraystretch{2}
\begin{table}
\caption{Probabilities in example 2\label{Ex3}}
\centering
\begin{tabular}{|cccccc|}
  \hline
   & \multicolumn{2}{c}{$P(S=1|U,T)$}  && \multicolumn{2}{c|}{$P(Y=1|U,S)$}   \\
\cline{2-3}\cline{5-6}
    & $T=0$ & $T=1$ && $S=0$ & $S=1$ \\
    \hline
  $U=0$ & 0.28 & 0.77 && 0.13 & 0.87 \\
  $U=1$ & 0.32 & 0.65 && 0.33 & 0.52 \\
  \hline
\end{tabular}
\end{table}
\renewcommand\arraystretch{0.8}
\vspace{-2mm}
\item[Example 3: (almost necessity)] Let $P(T=1)=0.5$, $P(U=1)=0.5$ and the conditional probabilities are given in Table \ref{Ex3_1}. Then we have $\mathrm{ACE}(T\rightarrow S)=0.10$, $\gamma=\mathrm{ACE}(S\rightarrow Y)=0.20$, $\mathrm{ACE}(T\rightarrow Y)=0.14$, and the right hand side in (\ref{SN}) is 0.60.  Thus the inequality (\ref{SN}) does not hold and the surrogate paradox does not occur. Note that we have the observed data probabilities: $P(T=1)=0.5$, $P(S=1|T=1)=0.7$,  $P(S=1|T=0)=0.6$, $P(Y=1,S=0|T=0)=0.1$, $P(Y=1,S=1|T=0)=0.3$ and $\gamma=0.2$. Now we find a data generating mechanism that yields the same observed data but suffers from the surrogate paradox. Consider $P(T=1)=0.5$, $P(U=1)=0.5$ and the conditional probabilities are given in Table \ref{Ex3_2}.  Then we have $\mathrm{ACE}(T\rightarrow S)=0.10$, $\gamma=\mathrm{ACE}(S\rightarrow Y)=0.20$, $\mathrm{ACE}(T\rightarrow Y)=-0.10$, and the right hand side in (\ref{SN}) is still 0.60.  Thus the inequality (\ref{SN}) does not hold and the surrogate paradox occurs.
\renewcommand\arraystretch{2}
\begin{table}
\caption{The first set of probabilities in example 3\label{Ex3_1}}
\centering
\begin{tabular}{|cccccc|}
  \hline
   & \multicolumn{2}{c}{$P(S=1|U,T)$}  && \multicolumn{2}{c|}{$P(Y=1|U,S)$}   \\
\cline{2-3}\cline{5-6}
    & $T=0$ & $T=1$ && $S=0$ & $S=1$ \\
    \hline
  $U=0$ & 0.70 & 0.50 && 0.60 & 0.40 \\
  $U=1$ & 0.50 & 0.90 && 0.40 & 0.64 \\
  \hline
\end{tabular}
\end{table}
\renewcommand\arraystretch{0.8}

\renewcommand\arraystretch{2}
\begin{table}
\caption{The second set of probabilities in example 3\label{Ex3_2}}
\centering
\begin{tabular}{|cccccc|}
  \hline
   & \multicolumn{2}{c}{$P(S=1|U,T)$}  && \multicolumn{2}{c|}{$P(Y=1|U,S)$}   \\
\cline{2-3}\cline{5-6}
    & $T=0$ & $T=1$ && $S=0$ & $S=1$ \\
    \hline
  $U=0$ & 0.70 & 0.50 && 0.20 & 0.80 \\
  $U=1$ & 0.50 & 0.90 && 0.28 & 0.08 \\
  \hline
\end{tabular}
\end{table}
\renewcommand\arraystretch{0.8}

\end{description}

\section{Extensions}\label{nonboundsection}
\subsection{Non-strong surrogate case}\label{nonstrong}
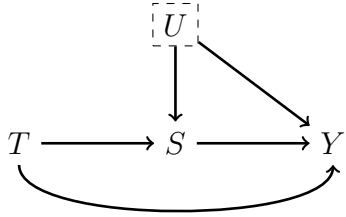
\begin{figure}
\centering
\begin{tikzpicture}
\node[text centered] (t) {$T$};
\node[right = 1.5 of t, text centered] (s) {$S$};
\node[right=1.5 of s, text centered] (y) {$Y$};
\node[draw, rectangle, dashed, above = 1 of s, text centered] (u) {$U$};

\draw[->, line width= 1] (t) --  (s);
\draw [->, line width= 1] (s) -- (y);
\draw[->,line width= 1] (u) --(s);
\draw[->,line width= 1] (u) -- (y);
\draw[->,line width=1] (t) to  [out=270,in=270, looseness=0.5] (y);
\end{tikzpicture}
\medskip
\caption{Causal diagram of non-strong surrogate  $S$  for the effect of the treatment $T$ on the outcome $Y$.\label{DAGnon}}
\end{figure}

At times, the treatment may have several causal pathways to the primary outcome that are not fully mediated through the surrogate endpoint. In that case, the assumption \ref{A2} does not hold (see Figure \ref{DAGnon}). For example, physical exercises decrease the cardiac output and ejection fraction in patients with severe heart failure but it may also affect the survival rate in other pathways such as lowering the risk of type 2 diabetes and some cancers.

When the assumption 2 is violated, the potential outcome $Y_{ts}$ can no longer be simplified as $Y_{s}$. There are 4 potential outcomes $Y_{ts}$ and 16 different possible values of the vector $(Y_{00},Y_{01},Y_{10},Y_{11})$ with $Y_{ts}\in \{0,1\}$ for $t=0,1$ and $s=0,1$. The number of the potential outcomes for surrogate remains at 2 and the different possible values for $(S_0,S_1)$ remains at 4 as in the strong surrogate case with $S_t\in\{0,1\}$ for $t=0,1$. Thus, we have $16\times 4=64$ different possible values for the vector $(Y_{00},Y_{01},Y_{10},Y_{11},S_{0},S_{1})$. Let $q_{i,j}\geq0$ denote the proportion of getting each possible value of $(Y_{00}, Y_{01}, Y_{10}, Y_{11}, S_{0},S_{1})$ in the whole population for $i=0,\ldots,15$ and $j=0,\ldots,3$, as shown in Table 1 in the supplementary material.

In the strong surrogate case, $\gamma=P(Y_{S=1}=1)-P(Y_{S=0}=1)$ describes the causal effect of the surrogate on the outcome. In contrast, in the non-strong surrogate case, we define the following quantities for the effects of the surrogate on the outcome:
$$\left\{
  \begin{array}{ll}
  \gamma_{0}=P(Y_{01}=1)-P(Y_{00}=1),\\
  \gamma_{1}=P(Y_{11}=1)-P(Y_{10}=1).   \\
  \end{array}
\right.$$
By definition, the $\gamma_{0}$ and $\gamma_{1}$ are the causal effects of $S$ on $Y$ when the treatment $T$ is set to 0 and 1 respectively. The point estimates or the ranges of $\gamma_{0}$ and $\gamma_{1}$ may be obtained from an external study or subject matter knowledge. 

Similar as in (\ref{constraint1}), we have the following linear constraints for the non-strong surrogate under the assumption \ref{A1}:
\begin{equation}\label{C2}
\left\{
  \begin{array}{ll}
    P(Y=0,S=0|T=0)=\sum_{i=0}^{7}(q_{i,0}+q_{i,1}),  \\
    P(Y=1,S=0|T=0)=\sum_{i=8}^{15}(q_{i,0}+q_{i,1}),  \\
    P(Y=0,S=1|T=0)=\sum_{i=0,1,2,3,8,9,10,11}(q_{i,2}+q_{i,3}),  \\
    P(S=0|T=1)=\sum_{i=0}^{15}(q_{i,0}+q_{i,2})  \\
  \gamma_{0}=\sum_{j=0,1,2,3}(q_{4,j}+q_{5,j}+q_{6,j}+q_{7,j}-q_{8,j}-q_{9,j}-q_{10,j}-q_{11,j})  ,  \\
  \gamma_{1}=\sum_{j=0,1,2,3}(q_{1,j}+q_{5,j}+q_{9,j}+q_{13,j}-q_{2,j}-q_{6,j}-q_{10,j}-q_{14,j}), \\
  1=\sum_{i=0}^{15}\sum_{j=0}^{3}q_{i,j}.
  \end{array}
\right.
\end{equation}
The parameter of interest $\mathrm{ACE}(T\rightarrow Y)$ can be expressed as a linear combination of $\{q_{i,j},i=0,\cdots,15,j=0,\cdots,3\}$ as:
\begin{eqnarray*}
  \mathrm{ACE}(T\rightarrow Y) &=& ~~P(Y_{T=1,S=0}=1,S_{T=1}=0)+P(Y_{T=1,S=1}=1,S_{T=1}=1) \\
   &&-P(Y_{T=0,S=0}=1,S_{T=0}=0)-P(Y_{T=0,S=1}=1,S_{T=0}=1)  \\
&=&~~ (q_{2,0}+q_{3,0}+q_{6,0}+q_{7,0}-q_{8,0}-q_{9,0}-q_{12,0}-q_{13,0})\\
&&+(q_{1,1}+q_{3,1}+q_{5,1}+q_{7,1}-q_{8,1}-q_{10,1}-q_{12,1}-q_{14,1})\\
&&+(q_{2,2}+q_{3,2}+q_{10,2}+q_{11,2}-q_{4,2}-q_{5,2}-q_{12,2}-q_{13,2})\\
&&+(q_{1,3}+q_{3,3}+q_{9,3}+q_{11,3}-q_{4,3}-q_{6,3}-q_{12,3}-q_{14,3}).
\end{eqnarray*}

Again, we use the linear programming method to obtain the sharp bounds for $\mathrm{ACE}(T\rightarrow Y)$ in the following theorem.

\vspace{-3mm}
\begin{Theorem}\label{Th4}
Under the assumption \ref{A1}, if $\big\{P(S,Y|T=0),P(S|T=1),\gamma_0,\gamma_1\big\}$ is known, we can have the sharp bounds for $\mathrm{ACE}(T\rightarrow Y)$, denote as $[L',U']$, where
$$L'= \max\left(
               \begin{array}{c}
                 -P(Y=1|T=0) \\
                 -\gamma_{1}-P(Y=1|T=0)-P(S=1|T=1) \\
                 \gamma_{1}-P(Y=1|T=0)-P(S=0|T=1) \\
               \end{array}
             \right)
$$
\noindent and

$$U'= \min\left(
               \begin{array}{c}
                 P(Y=0|T=0) \\
                 -\gamma_{1}+P(Y=0|T=0)+P(S=1|T=1) \\
                 \gamma_{1}+P(Y=0|T=0)+P(S=0|T=1) \\
               \end{array}
             \right).
$$
\end{Theorem}

Similar as in Theorem \ref{Th1}, the bounds in Theorem \ref{Th4} are also sharp in the sense that they can be obtained. Note that the bounds of $\mathrm{ACE}(T\rightarrow Y)$ do not involve $\gamma_{0}$. That is, if we delete the fifth constraint in (\ref{C2}), the bounds of $\mathrm{ACE}(T\rightarrow Y)$ will not change. To see this, note that due to the randomization of treatment $T$, we have $\mathrm{ACE}(T\rightarrow Y)=P(Y=1|T=1)-P(Y=1|T=0)$. Thus, bounding $\mathrm{ACE}(T\rightarrow Y)$ is equivalent to bounding $P(Y=1|T=1)$ since $P(Y=1|T=0)$ can be identified using the observed data. Also $P(Y=1|T=1)=P(Y_{10}=1,S_{1}=0)+P(Y_{11}=1,S_{1}=1)$, which only involves $(Y_{10},Y_{11},S_{1})$. However, $\gamma_{0}$ only involves the distributions of $Y_{00}$ and $Y_{01}$, which do not have any intersection with the parameters involved in $P(Y=1|T=1)$.

Also note that the bounds of $\mathrm{ACE}(T\rightarrow Y)$ only depend on the data from the control arm via $P(Y|T=0)$ rather than $P(Y,S|T=0)$. This is because in the non-strong surrogate case, bounding $P(Y=1|T=1)$ does not depend on $P(Y,S|T=0)$. The information of the controlled arm is only utilized when estimating $P(Y=1|T=0)$. That is, for the non-strong surrogate case, we only need to observe the primary endpoint $Y$ from the previous study ($T=0$) rather than jointly observe the surrogate $S$ and the outcome $Y$. The knowledge of the latter does not further shrink the bounds of $\mathrm{ACE}(T\rightarrow Y)$.

Since $\mathrm{ACE}(T\rightarrow Y)\geq L^{'}$, the condition that $L^{'}>0$ excludes the surrogate paradox. That is also equivalent to
\begin{equation}\label{non_SN}
\gamma_1>P(Y=1|T=0)+P(S=0|T=1).
\end{equation}
Thus, \eqref{non_SN} can be used as a criterion to exclude the paradox with a non-strong surrogate. We have the following theorem.

\begin{Theorem}\label{Th5}
Under the assumption \ref{A1}, if $\big\{P(Y|T=0),P(S|T=1),\gamma_1\big\}$ is known, the surrogate paradox can be excluded if (\ref{non_SN}) holds. Additionally, \eqref{non_SN} is an optimal criterion to exclude the surrogate paradox.
\end{Theorem}

As mentioned, when the surrogate endpoint is actually a strong surrogate, $Y_{st}=Y_s$, thus $\gamma_0=\gamma_1=\gamma$. Comparing $L^{'}$ with $L$, $[L',U']$ are wider than the bounds $[L,U]$ if replacing $\gamma_1$ with $\gamma$ in the bounds $L'$ and $U'$. Consequently, the right hand side of \eqref{non_SN} is greater than \eqref{SN}. That is, if the surrogate endpoint is actually a strong surrogate but we derived the bounds assuming it is a non-strong surrogate, then we would get more conservative results. 

\subsection{Causal relative risk}
Compared with the ACE which defines causal effect on the difference scale, the causal relative risk (CRR) defines causal effect on the relative risk scale. The CRR is popularly used in the statistical analysis of binary outcomes where the outcome of interest has a relatively low probability, such as having some rare disease.

For the simplicity of illustration, we only discuss the extension among the class of strong surrogate endpoints. The extension for non-strong surrogate can be generalized in a similar fashion as in the ACE scale. Under the assumption 2, let $\mathrm{CRR}(T\rightarrow Y)=P(Y_{T=1}=1)/P(Y_{T=0}=1)$ denote the CRR of treatment $T$ on the outcome $Y$. Similarly define $\mathrm{CRR}(T\rightarrow S)$ and $\mathrm{CRR}(S\rightarrow Y)$.

Note that $\mathrm{CRR}(T \rightarrow Y)>1$ is equivalent to $P(Y_{T=1}=1)>P(Y_{T=0}=1)$, which indicates a positive treatment effect on the outcome. Similarly, $\mathrm{CRR}(T \rightarrow Y)<1$ indicates a negative treatment effect on the outcome. Thus, we can restate the phenomenon of surrogate paradox on the CRR scale as:

\begin{itemize}
  \item A treatment has a positive effect on a surrogate on the CRR scale which in turn has a positive effect on a primary outcome on the CRR scale, but the treatment has a negative effect on the primary outcome on the CRR scale, i.e., $\mathrm{CRR}(T\rightarrow S)>1$, and $\mathrm{CRR}(S\rightarrow Y)>1$ but $\mathrm{CRR}(T\rightarrow Y)<1$.
\end{itemize}

Similar to the case of ACE, we also have other surrogate paradox phenomenon on the CRR scale by redefining $T^{*}=1-T$, $S^{*}=1-S$, $Y^{*}=1-Y$, which we omitted here.


To avoid the surrogate paradox, we first bound $\mathrm{CRR}(T\rightarrow Y)$ with the available data $\{P(Y,S|T=0),P(S|T=1),\mathrm{CRR}(S\rightarrow Y)\}$. With the same notation of potential outcomes types in the ACE case and let $\gamma_{CRR}=\mathrm{CRR}(S\rightarrow Y)$, we can have the bounding problem under the assumptions \ref{A1} and \ref{A2}:
$$\mathrm{bound}~\mathrm{CRR}(T\rightarrow Y)=\frac{q_{02}+q_{03}+q_{22}+q_{23}+q_{11}+q_{13}+q_{31}+q_{33}}{q_{02}+q_{03}+q_{12}+q_{13}+q_{21}+q_{23}+q_{31}+q_{33}},$$
subjects to
\begin{equation}\label{constraint2}
\left\{
  \begin{array}{ll}
    P(Y=0,S=0|T=0)=q_{00}+q_{01}+q_{10}+q_{11},  \\
    P(Y=1,S=0|T=0)=q_{02}+q_{03}+q_{12}+q_{13},  \\
        P(Y=0,S=1|T=0)=q_{20}+q_{22}+q_{30}+q_{32},  \\
    P(S=0|T=1)=q_{00}+q_{01}+q_{20}+q_{21}+q_{02}+q_{03}+q_{22}+q_{23},  \\
    \begin{split}\gamma_{CRR}=\frac{q_{01}+q_{11}+q_{21}+q_{31}+q_{03}+q_{13}+q_{23}+q_{33}}{q_{02}+q_{12}+q_{22}+q_{32}+q_{03}+q_{13}+q_{23}+q_{33}},
    \end{split} \\
    \sum_{i=0}^{3}\sum_{j=0}^{3}q_{ij}=1, \\
    q_{ij}\geq0,i,j=0,1,2,3.
  \end{array}
\right.
\end{equation}
Note that both the objective function and the constraints are nonlinear, thus the linear programming cannot be applied directly. However, note that the fifth constraint $\gamma_{CRR}=\big(q_{01}+q_{11}+q_{21}+q_{31}+q_{03}+q_{13}+q_{23}+q_{33}\big)/\big(q_{02}+q_{12}+q_{22}+q_{32}+q_{03}+q_{13}+q_{23}+q_{33}\big)$ in (\ref{constraint2}) can be expressed as a linear constraint:
$$\big(q_{01}+q_{11}+q_{21}+q_{31}+q_{03}+q_{13}+q_{23}+q_{33}\big)-\big(q_{02}+q_{12}+q_{22}+q_{32}+q_{03}+q_{13}+q_{23}+q_{33}\big)\gamma_{CRR}=0.$$
In addition, the objective function is still not linear. The following theorem converts a linear fractional programming problem into a linear programming problem:

\begin{Theorem}\label{Th6}
\cite{schaible1983fractional}. {\it The linear fractional programming problem
\begin{equation}\label{citeTh}
\sup\{\frac{\gamma^{T}x+a}{\alpha^{T}x+b}:~Ax\leq \beta, \alpha^{T}x+b>0,x\geq0\}
\end{equation}
is equivalent to the linear programming problem
\begin{equation*}
\sup\{\gamma^{T}y+at:~Ay-\beta t\leq0,\alpha^{T}y+bt=1,y\geq0,t>0\}
\end{equation*}
by setting $y=x/(\alpha^{T}x+b)$ and $t=1/(\alpha^{T}x+b)$, where $A$ is a coefficient matrix and $\alpha,\beta,\gamma$ are constant vectors; $x$ is a vector of variables; $a$ and $b$ are constants; and $x\geq0$ indicates that each component of $x$ is larger or equal to zero. The strict inequality $t>0$ can be replaced by $t\geq0$ if equation (\ref{citeTh}) has an optimal solution.}
\end{Theorem}


The criterion to exclude the surrogate paradox in the relative risk scale can be derived by setting the lower bound to be greater than 1. Furthermore, such criterion also satisfies the optimality given in definition \ref{def: optimal}.

\section{Statistical Analysis}

We now apply the methods developed in Sections \ref{boundsection} and \ref{nonboundsection} to investigate the use of hypertension as a valid surrogate to evaluate the effect of anti-hypertension drugs on the mortality. In medicine, systolic hypertension is defined as an elevated SBP $>$140 mm Hg and our surrogate is defined as a binary variable for the hypertension condition ($S=1$ for normotensive and $S=0$ for hypertensive). As reported \cite{accord2010effects}, the average SBP was stabilized after 1 year for both intensive and standard therapy groups. Since later visits suffer from lower attenuation rate, we choose the SBP at 1 year as the surrogate for the ACCORD study. The outcome $Y$ is also a binary variable indicating death at 10 years ($Y=1$ for survived and $Y=0$ for death). Under the potential outcomes framework, Lawes et al. \cite{lawes2008global} reported 13.5\% of all deaths was due to high blood pressure. Since Lawes et al.  \cite{lawes2008global} did not report the percentage reduction of death due to high blood pressure among diabetic patients, as an illustration, we assume for the diabetic subgroup, the percentage reduction was also 13.5\%. This information could be more accurately collected in subgroup analysis.


Since both therapies were designed to lower the blood pressure to certain level, we first assume that the SBP is a strong surrogate. Assuming the rate of death among normotensive (number of death divided by number of normotensive people) is $P(Y_{S=0}=1)=0.1$, we obtain that on the difference scale, the causal effect of $S$ on $Y$ is $\gamma=0.0135$. By Theorem \ref{Th1}, we obtained the bounds for $\text{ACE}(T\xrightarrow{}Y)$ are [-0.331, 0.165]. Using the bootstrap method, we obtained the estimated standard deviation for upper bound and lower bound to be 0.009 and 0.006 respectively, thus the uncertainty region (confidence interval for bounds, see more in \cite{richardson2014nonparametric}) is [-0.349, 0.173]. Thus, both the bounds and the uncertainty region includes 0, that is the surrogate paradox cannot be excluded without further assumptions. Note that in the calculation of the uncertainty region, we treated $\text{ACE}(S\xrightarrow{}Y)=\gamma$ as a known number since the information is collected from a global investigation \cite{lawes2008global}. If we obtained $\gamma$ from a finite sample study, we could apply bootstrap method to that study to incorporate the variability of the estimate of $\gamma$.

We also utilized other value for $P(Y_{S=0}=1)$ but they all resulted in bounds of ACE($T\xrightarrow{}Y$) containing 0, indicating the sign of $\text{ACE}(T\xrightarrow{}Y)$ depends on the joint distribution of unobserved confounder $U$ and observed variables. Thus, with the observed data, for the bounds of ACE($T\xrightarrow{}Y$) to exclude 0, one needs to have $\gamma>0.688$. 

We then relax the strong surrogate assumption to derive bounds for $\text{ACE}(T\xrightarrow{}Y)$. As indicated in Section \ref{nonstrong}, the information of $\gamma_1$ can be obtained by an external study or subject matter knowledge. 
For illustration purpose, we first assume $\gamma_1=\gamma=0.0135$. Thus the bounds for $\text{ACE}(T\xrightarrow{}Y)$ are [-0.835, 0.165] with bootstrap uncertainty region [-0.846, 0.175]. As expected, allowing the surrogate endpoint to be non-strong when $\gamma_1=\gamma$ yields wider bounds. From the lower bound $L'$, we need $\gamma_{1}> 0.927$ to yield lower bound greater than 0. Such large $\gamma$ is impossible to reach for a non-pandemic disease.

Finally, we derive bounds for $\text{CRR}(T\xrightarrow{}Y)$ when the contrast is on the causal relative risk scale and evaluate it using the R package ``linprog" \cite{linprog2012}. Again, we make the strong surrogate assumption. Note that Lawes et al. \cite{lawes2008global} indicated that $\gamma_{\mathrm{CRR}}=1.135$. Thus, the bounds for $\text{CRR}(T\xrightarrow{}Y)$ are [0.604, 1.197] with bootstrap uncertainty region [0.583, 1.210]. The bounds include 1 which also indicates the surrogate paradox cannot be avoided in this case. 

Thus, we conclude that for evaluating the effect of anti-hypertension drug on the long-term death, using high blood pressure as a surrogate cannot guarantee the bounds to exclude null. That is, if the unmeasured confounders have certain value, it is possible that the treatment has a possible effect in reducing the high blood pressure and lowering the high blood pressure could reduce the death rate, but the treatment could increase the death rate. Thus, for the development of such anti-hypertension drug, it is recommended to also collect the information on the long-term death rate.

\section{Discussion}
In this article, we derived sharp bounds for the ACE of treatment on the outcome with a surrogate endpoint. Based on the sharp bounds, we proposed novel criteria to exclude surrogate paradox. Note that both our bounds and the conditions to avoid surrogate paradox are based on the observed data as opposed to the full data $(U,T,S,Y)$. This allows us to exclude the surrogate paradox only using the observed data.



In the strong surrogate scenario, the information of $\gamma$ is important for the results derived. When surrogate endpoint can be manipulated directly, measuring it would be possible. When the surrogate is not that easy to manipulate, various methods could be used to estimate or obtain the bounds of $\mathrm{ACE}(S\rightarrow Y)$ in the presence of unmeasured confounding. Once the information of $\gamma$ is collected, it can be used later to evaluate the effect of different treatment on the same surrogate and primary outcome. Such transportability is appealing since it may save the time and cost of future clinical trials.

If the strong surrogate assumption is violated, we need the information for $\gamma_{1}$, which may be more challenging to collect than $\gamma$, but a sensitivity analysis could be carried out if we know a range of $\gamma_1$ based on the subject matter knowledge. Alternatively, one could try to find a vector of surrogates so that all the pathways of the effect of treatment to the primary outcome is mediated through the vector of surrogates. Then the strong surrogate assumption may hold for this vector of surrogates. The composition of such surrogate vector depends on the subject matter knowledge of the underlying data generating mechanism. However, as the dimension of the surrogates vector increases, the computation burden increases exponentially, which is left for future discovery.

Similar methods could be modified for multinomial outcome $Y$, however, computational challenge also increases exponentially. We also leave the conditions to avoid surrogate paradox for continuous $Y$ and survival outcome for future research topic.

In this article, we mainly focused on the exclusion of the surrogate paradox when the surrogate fully mediates or partially mediates the effect of the treatment on the outcome of interest. Vanderweele \cite{vanderweele2013surrogate} considered extending the surrogate criteria when the surrogate is not on the causal pathway from the treatment to the outcome. He suggested that the confounding between the surrogate and the outcome may provide a source of importance information of surrogacy. To the extreme, he suggested that even the surrogate has no causal effect on the outcome, a surrogate can still be a good surrogate (See figure 3 in \cite{vanderweele2013surrogate}). Also, in the definition of surrogate paradox, he replaced the condition of $\mathrm{ACE}(S\rightarrow Y)>0$ with $S$ and $Y$ are positively associated. We leave extending our criteria to avoid unmediated surrogate paradox as future research topic.


\section*{Acknowledge}
We thank Haoda Fu, Linbo Wang and Wei Sun for helpful discussion.

\bibliographystyle{plain}

\bibliography{Arxiv_surrogate}

\end{document}


\maketitle

\markboth{\hfill{\footnotesize\rm YIN ET AL.} \hfill}
{\hfill {\footnotesize\rm NOVEL CRITERIA TO EXCLUDE SURROGATE PARADOX} \hfill}


\setcounter{section}{0}
\setcounter{equation}{0}
\def\theequation{S\arabic{section}.\arabic{equation}}
\def\thesection{S\arabic{section}}

\section{Proof of Theorem 1}
Let $\overrightarrow{q}=(q_{00},q_{01},q_{02},q_{03},q_{10},q_{11},q_{12},q_{13},q_{20},q_{21},q_{22},q_{23},q_{30},q_{31},q_{32},q_{33})'$,
$b=(P(Y=0,S=0|T=0),P(Y=1,S=0|T=0),P(Y=0,S=1|T=0),P(S=0|T=1),\gamma,1)'$, where the superscript $'$ denote the matrix transpose. Let
\begin{center}
$A=\left(
  \begin{array}{cccccccccccccccc}
     1&  1& 0 & 0 & 1 & 1 & 0 & 0 & 0 & 0 & 0 & 0 & 0 & 0 & 0 & 0 \\
     0& 0 & 1 &1  & 0 & 0 & 1 & 1 & 0 & 0 & 0 & 0 & 0 & 0 & 0 & 0 \\
     0& 0 & 0 & 0 & 0 & 0 & 0 & 0 & 1 &0  & 1 & 0 & 1 & 0 & 1 & 0 \\
     1& 1 & 1 & 1 & 0 & 0 & 0 & 0 & 1 & 1 & 1 & 1 & 0 & 0 & 0 & 0 \\
     0& 1 &-1 & 0 & 0 & 1 &-1 & 0 & 0 & 1 &-1 & 0 & 0 & 1 &-1 & 0 \\
     1& 1 & 1 & 1 & 1 & 1 & 1 & 1 & 1 & 1 & 1 & 1 & 1 & 1 & 1 & 1 \\
  \end{array}
\right),$
\end{center}
then all the constraints in (1) can be written as $A\overrightarrow{q}=b$. Note that $\mathrm{ACE}(T\rightarrow Y)=c'\overrightarrow{q}$, where $c=(0, 0, 0, 0, 0,1,-1,0,0,-1,1,0,0,0,0,0)'$. Thus the problem of obtaining the lower bound for $\mathrm{ACE}(T\rightarrow Y)$ can be equivalently rewritten as a minimization problem with linear constraints:
\begin{equation*}\label{LP}
\min~ c'\overrightarrow{q},~~\mathrm{subject~ to}~~ \left\{
                                            \begin{array}{ll}
                                             A\overrightarrow{q}=b ,  \\
                                             \overrightarrow{q}\geq0 .
                                            \end{array}
                                          \right.
\end{equation*}
The upper bound can be obtained similarly by replacing min with max in the objective function and then use the dual linear programming \citep{dantzig2003linear}. 

To illustrate the derivation of the bounds, we consider a general linear programming problem with the form:
\begin{flushright}
bound $c'\overrightarrow{q}$,~~subject to$~~ \left\{
                                            \begin{array}{ll}
                                             A\overrightarrow{q}=b ,  \\
                                             \overrightarrow{q}\geq0 ,
                                            \end{array}
                                          \right.$~~~~~~~~~~~~~~~~~~~~~~~~~~~~~~~~~~~(P)
\end{flushright}
where $A$ is a $m\times n$ $(m<n)$ numeric matrix, $b$ is a $m\times 1$ symbolic vector, $c$ is a $n\times1$ numeric vector, and $\overrightarrow{q}$ is a $n\times 1$ vector.
%

Note the problem (P) can be divided into two problems:
\begin{flushright}
min $c'\overrightarrow{q}$,~~subject to$~~ \left\{
                                            \begin{array}{ll}
                                             A\overrightarrow{q}=b ,  \\
                                             \overrightarrow{q}\geq0 ,
                                            \end{array}
                                          \right.$~~~~~~~~~~~~~~~~~~~~~~~~~~~~~~~~~~~(P1)
\end{flushright}
and
\begin{flushright}
max $c'\overrightarrow{q}$,~~subject to$~~ \left\{
                                            \begin{array}{ll}
                                             A\overrightarrow{q}=b ,  \\
                                             \overrightarrow{q}\geq0 ,
                                            \end{array}
                                          \right.$~~~~~~~~~~~~~~~~~~~~~~~~~~~~~~~~~~~(P2)
\end{flushright}
To solve the minimization problem (P1), we first focus on solving its dual problem (D1), where (D1) is

\vspace{-2mm}\begin{flushright}
max $b'\overrightarrow{p}$,~~subject to~~$  A'\overrightarrow{p}\leq c .$~~~~~~~~~~~~~~~~~~~~~~~~~~~~~~~~~~~(D1)
\end{flushright}

\vspace{-2mm}The set $Q=\{\overrightarrow{p}|A'\overrightarrow{p}\leq c\}$ is polyhedron and the quantity $b'\overrightarrow{p}$ reaches its extreme value at the vertexes of the polyhedron $Q$. Since both $A$ and $c$ do not involve any symbolics, we can enumerate all the vertexes $\{p_{1},\cdots,p_{K}\}$ of $Q$. Thus, the solution to the optimal problem (D1) is
$$L=\max\{b'p_{1},\cdots,b'p_{K}\}.$$
The bound $L$ can be attained since all the vertexes belong to the set $Q$.

According to \cite{dantzig2003linear}, an optimal solution $b'\overrightarrow{p*}$ to the dual problem (D1) corresponds to an optimal solution $c'\overrightarrow{q*}$ to the primal problem (P1) with $c'\overrightarrow{q*}=b'\overrightarrow{p*}$. Thus, we can conclude that the optimal solution to (P1) is $L=\max\{b'p_{1},\cdots,b'p_{K}\},$ and this can be attained.

Similarly, we can obtain the optimal solution for the problem (P2), denote as $U$, then the optimal problem to (P) is
$$L\leq c'\overrightarrow{q}\leq U,$$
and the bounds are sharp.

%

\vspace{-7mm}\section{Proof of Theorem 2 and 3}
The proof of theorem 2 follows from the fact that $L>0$ is a sufficient condition for $\mathrm{ACE}(T\rightarrow Y)>0$. The optimality of condition (2) follows directly from the sharpness of the bounds given in theorem 1.

\section{Examination of Optimality for the Existing Criteria}
To our knowledge, the criteria to preclude the surrogate paradox has been considered by \citet{chen2007criteria}, \citet{ju2010criteria}, \citet{wu2011sufficient} and \citet{vanderweele2013surrogate}, where the former two were also summarized in \citet{vanderweele2013surrogate}. We show that their criteria are sufficient but none is optimal to exclude the surrogate paradox.

\citet{chen2007criteria} assumed that if the following two conditions hold for a strong surrogate $S$ then the surrogate paradox can be excluded: (a) $E(Y|S,U)$ monotonic in $S$ and (b) $P(S>s|t,u)$ is increasing in $t$ for all $s$ and $u$. Suppose $P(U=1)=0.5$ and the rest probabilities are given in table \ref{ExS1}. It can be verified that $\gamma=0.8>0.625$, where the latter is the right hand side of inequality (2). Thus, by theorem 2, there exists no underlying full data that has surrogate paradox. Note that this example also violates Chen's condition (b). However, for all full data $(U,T,S,Y)$ with the same observed data, none of them has the surrogate paradox thus violates (ii) of the optimality definition 1.

When the outcome is binary, the conditions of \citet{ju2010criteria} are equivalent to the conditions of \citet{chen2007criteria}. Thus, the counter example mentioned above for \citet{chen2007criteria} also serves as a counter example for \citet{ju2010criteria}. \citet{vanderweele2013surrogate} generalizes the results of \citet{ju2010criteria} to allow for nonstrong surrogate. Note that \citet{vanderweele2013surrogate} also assumes $P(S>s|t,u)$ is increasing in $t$ for all $s$ and $u$. Thus, this example also violates Vanderweele's condition.

We next examine the optimality of \citet{wu2011sufficient}. \citet{wu2011sufficient} assumed that if following two conditions hold, then the surrogate paradox can be excluded: (a) $E(Y|S=s,T=1)$ or $E(Y|S=s,T=0)$ monotonically increases as $s$ increases; and (b) $E(Y|S=s,T=1)>E(Y|S=s,T=0)$ for all $s$. Assume $P(U=1)=0.3$ and the rest probabilities are given in table \ref{ExS2}. Note that $E(Y|S=1,T=1)=0.828$ and $E(Y|S=1,T=0)=0.834$ thus the condition (b) of Wu does not satisfy. However, $\gamma=0.7>0.572$, where the latter is the right hand side of inequality (2). Thus, by theorem 2, there exists no underlying full data that has surrogate paradox indicating (ii) of the optimality definition 1 is violated.

\renewcommand\arraystretch{2}
\begin{table}
\caption{Counter example for optimality of the criteria proposed by \citet{chen2007criteria}}\label{ExS1}
\begin{center}
\begin{tabular}{|cccccc|}
  \hline
   & \multicolumn{2}{c}{$P(S=1|U,T)$}  && \multicolumn{2}{c|}{$P(Y=1|U,S)$}   \\
\cline{2-3}\cline{5-6}
    & $T=0$ & $T=1$ && $S=0$ & $S=1$ \\
    \hline
  $U=0$ & 0.40 & 0.30 && 0.10 & 0.90 \\
  $U=1$ & 0.10 & 0.90 && 0.10 & 0.90 \\
  \hline
\end{tabular}
\end{center}
\end{table}
\renewcommand\arraystretch{0.8}

\vspace{-7mm}\renewcommand\arraystretch{2}
\begin{table}
\caption{Counter example for optimality of the criteria proposed by \citet{wu2011sufficient}}\label{ExS2}
\begin{center}
\begin{tabular}{|cccccc|}
  \hline
   & \multicolumn{2}{c}{$P(S=1|U,T)$}  && \multicolumn{2}{c|}{$P(Y=1|U,S)$}   \\
\cline{2-3}\cline{5-6}
    & $T=0$ & $T=1$ && $S=0$ & $S=1$ \\
    \hline
  $U=0$ & 0.50 & 0.90 && 0.10 & 0.80 \\
  $U=1$ & 0.60 & 0.80 && 0.20 & 0.90 \\
  \hline
\end{tabular}
\end{center}
\end{table}
\renewcommand\arraystretch{0.8}

\section{Bounds Without the Magnitude of $\gamma$}
If we only have the available data on $\{P(Y,S|T=0), P(S|T), a<\gamma<b\}$, since the bounds we derived in theorem 1 are sharp, we can have the following condition to exclude paradox:

\begin{equation}\label{SN_ab}
a>\min\left(
               \begin{array}{c}
                2P(Y=1,S=1|T=0)+P(Y=0,S=0|T=0)\\
                P(S=0|T=1)+P(Y=1,S=1|T=0)\\
               \end{array}
             \right).
\end{equation}

\noindent Similar as theorem 3, this also satisfies the optimality defined in definition 1. As a special example, when $a=0$, \eqref{SN_ab} cannot be satisfied regardless of the observed data. This can also be seen by deriving the sharp bounds for $\mathrm{ACE}(T\rightarrow Y)$ directly with the available data $\{P(Y,S|T=0), P(S|T), \gamma>0\}$:

\vspace{-7mm}$$\mathrm{ACE}(T\rightarrow Y)\geq \max\left(
               \begin{array}{c}
                 -P(Y=1,S=0|T=0)-P(S=0|T=1) \\
                 -P(Y=1|T=0) \\
                 -P(Y=1,S=1|T=0)-P(S=1|T=1) \\
                 -2P(Y=1,S=1|T=0)-P(Y=0,S=0|T=0) \\
                 -P(Y=1,S=1|T=0)-P(S=0|T=1) \\
               \end{array}
             \right),
$$
\noindent and
$$\mathrm{ACE}(T\rightarrow Y)\leq \min\left(
               \begin{array}{c}
                 P(Y=0,S=0|T=0)+P(S=0|T=1) \\
                 P(Y=0|T=0) \\
                 P(Y=0,S=1|T=0)+P(S=1|T=1) \\
                 2P(Y=0,S=0|T=0)+P(Y=1,S=1|T=0) \\
                 P(Y=0,S=0|T=0)+P(S=1|T=1) \\
               \end{array}
             \right).
$$
It is easy to see that the bounds above always contain zero.

\section{The Potential Outcomes Types for Non-strong Surrogate Scenario}
With a non-strong surrogate, we have 64 potential outcomes types in Table \ref{type} where $W=(Y_{00},Y_{01},Y_{10},Y_{11})$:
\renewcommand\arraystretch{1.6}
\begin{table}
\caption{The probabilities of the potential outcomes types for non-strong surrogate scenario where $W=(Y_{00},Y_{01},Y_{10},Y_{11})$}\label{type}
\begin{center}
\begin{tabular}{|l|l|l|l|l|}
  \hline
   & $S_{T=0}=0$ & $S_{T=0}=0$ & $S_{T=0}=1$ & $S_{T=0}=1$ \\
   & $S_{T=1}=0$ & $S_{T=1}=1$ & $S_{T=1}=0$ & $S_{T=1}=1$\\
    \hline
  $W=(0,0,0,0)$ & $q_{0,0}$ & $q_{0,1}$ & $q_{0,2}$ & $q_{0,3}$ \\
   \hline
  $W=(0,0,0,1)$ & $q_{1,0}$ & $q_{1,1}$ & $q_{1,2}$ & $q_{1,3}$ \\
   \hline
  $W=(0,0,1,0)$ & $q_{2,0}$ & $q_{2,1}$ & $q_{2,2}$ & $q_{2,3}$ \\
   \hline
  $W=(0,0,1,1)$ & $q_{3,0}$ & $q_{3,1}$ & $q_{3,2}$ & $q_{3,3}$ \\
  \hline
  $W=(0,1,0,0)$ & $q_{4,0}$ & $q_{4,1}$ & $q_{4,2}$ & $q_{4,3}$ \\
   \hline
  $W=(0,1,0,1)$ & $q_{5,0}$ & $q_{5,1}$ & $q_{5,2}$ & $q_{5,3}$ \\
   \hline
  $W=(0,1,1,0)$ & $q_{6,0}$ & $q_{6,1}$ & $q_{6,2}$ & $q_{6,3}$ \\
   \hline
  $W=(0,1,1,1)$ & $q_{7,0}$ & $q_{7,1}$ & $q_{7,2}$ & $q_{7,3}$ \\
  \hline
  $W=(1,0,0,0)$ & $q_{8,0}$ & $q_{8,1}$ & $q_{8,2}$ & $q_{8,3}$ \\
   \hline
  $W=(1,0,0,1)$ & $q_{9,0}$ & $q_{9,1}$ & $q_{9,2}$ & $q_{9,3}$ \\
   \hline
  $W=(1,0,1,0)$ & $q_{10,0}$ & $q_{10,1}$ & $q_{10,2}$ & $q_{10,3}$ \\
   \hline
  $W=(1,0,1,1)$ & $q_{11,0}$ & $q_{11,1}$ & $q_{11,2}$ & $q_{11,3}$ \\
  \hline
  $W=(1,1,0,0)$ & $q_{12,0}$ & $q_{12,1}$ & $q_{12,2}$ & $q_{12,3}$ \\
   \hline
  $W=(1,1,0,1)$ & $q_{13,0}$ & $q_{13,1}$ & $q_{13,2}$ & $q_{13,3}$ \\
   \hline
  $W=(1,1,1,0)$ & $q_{14,0}$ & $q_{14,1}$ & $q_{14,2}$ & $q_{14,3}$ \\
   \hline
  $W=(1,1,1,1)$ & $q_{15,0}$ & $q_{15,1}$ & $q_{15,2}$ & $q_{15,3}$ \\
  \hline
\end{tabular}
\end{center}
\end{table}

\bibliographystyle{apa}
\bibliography{/Users/Lan/Sync/Lan/study/statistics/biostatistics/research/my_paper/Surrogate_bounds/surrogate_manuscript03/surrogate}